\documentclass{article}\begin{document}\centerline{\bf\ Parametric Evaluations of the Rogers Ramanujan Continued Fraction}\vskip .10in

\centerline{\bf Nikos Bagis}

\centerline{Department of Informatics}

\centerline{Aristotele University of Thessaloniki}
\centerline{Thessaloniki, Greece}
\centerline{nikosbagis@hotmail.gr}

\begin{quote}
\begin{abstract}
In this article with the help of the inverse function of the singular moduli we evaluate the Rogers Ranmanujan continued fraction and his first derivative.  
\end{abstract}

\end{quote}

\label{intro}
\section{Introductory Definitions and Formulas}

For $\left|q\right|<1$, the Rogers Ramanujan continued fraction (RRCF) (see [6]) is defined as
\begin{equation}
R(q):=\frac{q^{1/5}}{1+}\frac{q^1}{1+}\frac{q^2}{1+}\frac{q^3}{1+}\cdots  
\end{equation}
We also define
\begin{equation}
(a;q)_n:=\prod^{n-1}_{k=0}{(1-aq^k)}
\end{equation}
\begin{equation}
f(-q):=\prod^{\infty}_{n=1}(1-q^n)=(q;q)_{\infty}
\end{equation}
Ramanujan give the following relations which are very useful:
\begin{equation}
\frac{1}{R(q)}-1-R(q)=\frac{f(-q^{1/5})}{q^{1/5}f(-q^5)} 
\end{equation} 
\begin{equation}
\frac{1}{R^5(q)}-11-R^5(q)=\frac{f^6(-q)}{q f^6(-q^5)} 
\end{equation} 
From the Theory of Elliptic Functions (see [6],[7],[10]),
\begin{equation}
K(x)=\int^{\pi/2}_{0}\frac{1}{\sqrt{1-x^2\sin(t)^2}}dt
\end{equation}
is the Elliptic integral of the first kind.
It is known that the inverse elliptic nome $k=k_r$, $k'^2_r=1-k^2_r$ is the solution of the equation
\begin{equation}
\frac{K\left(k'_r\right)}{K(k)}=\sqrt{r}
\end{equation}
where $r \in \bf R^{*}_+ \rm$. When $r$ is rational then the $k_r$ are algebraic numbers.\\
We can also write the function $f$ using elliptic functions. It holds (see [10]):
\begin{equation}
f(-q)^8=\frac{2^{8/3}}{\pi^4}q^{-1/3}(k_r)^{2/3}(k'_r)^{8/3}K(k_r)^4
\end{equation}
also holds 
\begin{equation}
f(-q^2)^6=\frac{2k_r k'_r K(k_r)^3}{\pi^3 q^{1/2}}
\end{equation}
From [5] it is known that
\begin{equation}
R'(q)=1/5q^{-5/6}f(-q)^4R(q)\sqrt[6]{R(q)^{-5}-11-R(q)^5}
\end{equation}
Consider now for every $0<x<1$ the equation $$x=k_r ,$$ which, have solution 
\begin{equation}
r=k^{(-1)}(x) 
\end{equation}
Hence for example $$k^{(-1)}\left(\frac{1}{\sqrt{2}}\right)=1$$ 
With the help of $k^{(-1)}$ function we evaluate the Rogers Ramanujan continued fraction.

\section{Propositions}

The relation between $k_{25r}$ and $k_r$ is (see [6] pg. 280):
\begin{equation}
k_rk_{25r}+k'_rk'_{25r}+2\cdot4^{1/3} (k_rk_{25r}k'_rk'_{25r})^{1/3}=1
\end{equation}
For to solve equation (12) we give the following
\[
\]  
\textbf{Proposition 1.}\\
The solution of the equation
\begin{equation}
x^6+x^3(-16+10x^2)w+15x^4w^2-20x^3w^3+15x^2w^4+x(10-16x^2)w^5+w^6=0
\end{equation}
when we know the $w$ is given by
$$
\frac{y^{1/2}}{w^{1/2}}=\frac{w^{1/2}}{x^{1/2}}=$$
\begin{equation}
=\frac{1}{2}\sqrt{4+\frac{2}{3}\left(\frac{L^{1/6}}{M^{1/6}}-4\frac{M^{1/6}}{L^{1/6}}\right)^2}+\frac{1}{2}\sqrt{\frac{2}{3}}\left(\frac{L^{1/6}}{M^{1/6}}-4\frac{M^{1/6}}{L^{1/6}}\right)
\end{equation}
where
\begin{equation}
w=\sqrt{\frac{L(18+L)}{6(64+3L)}}<1
\end{equation}
and
\begin{equation} 
M=\frac{18+L}{64+3L}
\end{equation}
If happens $x=k_r$ and $y=k_{25r}$, then $r=k^{(-1)}(x)$ and $w^2=k_{25r}k_r$, $(w')^2=k'_{25r}k'_r$.\\
\textbf{Proof.}\\
The relation (14) can be found using Mathematica. See also [6].
\[
\]
\textbf{Proposition 2.}\\
If $q=e^{-\pi\sqrt{r}}$ and 
\begin{equation}
a=a_r=\left(\frac{k'_r}{k'_{25r}}\right)^2\sqrt{\frac{k_r}{k_{25r}}}M_5(r)^{-3}
\end{equation}
Then 
\begin{equation}
a_r=R(q)^{-5}-11-R^5(q)
\end{equation}
Where $M_5(r)$ is root of: $(5x-1)^5(1-x)=256 (k_r)^2 (k'_r)^2 x$.\\
\textbf{Proof.}\\
Suppose that $N=n^2\mu$, where $n$ is positive integer and $\mu$ is positive real then it holds that
\begin{equation}
K[n^2\mu]=M_n(\mu)K[\mu]
\end{equation}
Where $K[\mu]=K(k_{\mu})$\\ 
The following formula for $M_5(r)$ is known
\begin{equation}
(5M_5(r)-1)^5(1-M_5(r))=256(k_r)^2 (k'_r)^2M_5(r)
\end{equation}
Thus if we use (5) and (8) and the above consequence of the Theory of Elliptic Functions, we get:
\[
R^{-5}(q)-11-R^{5}(q)=\frac{f^6(-q)}{q f^6(-q^5)}=a=a_r
\]
See also [4],[5].
\section{The Main Theorem}
From Proposition 2 and relation $w^2=k_{25r}k_r$ we get
\begin{equation}
w^5-k^2_rw=\frac{k^3_r(k^2_r-1)}{a_rM_5(r)^3} 
\end{equation}
Combining (13) and (21), we get: 
$$
[-10k_r^4+26k_r^6+a_rM_5(r)^3k_r^6-16k_r^8]+[-k_r^3-6a_rM_5(r)^3k_r^3+k_r^5-6a_rM_5(r)^3k_r^5]w+$$
\begin{equation}
+[a_rM_5(r)^3k_r^2+15a_rM_5(r)k_r^4]w^2-20a_rM_5(r)^3k_r^3w^3+15a_rM_5(r)^3k_r^2w^4=0
\end{equation}
Solving with respect to $a_rM_5(r)^3$, we get
\begin{equation}
a_r M_5(r)^3=\frac{16k_r^6-26k_r^4-wk_r^3+10k_r^2+wk_r}{k_r^4-6k_r^3w-20k_r^3w+15w^2k_r^2-6k_rw+15w^4+w^2}
\end{equation}
Also we have $$\frac{K(k_{25r})}{K(k_{r})}=M_5(r)=\frac{1}{m}=\left(\sqrt{\frac{k_{25 r}}{k_r}}+\sqrt{\frac{k'_{25r}}{k'_{r}}}-\sqrt{\frac{k_{25r}k'_{25r}}{k_rk'_r}}\right)^{-1}$$
$$=\left(\frac{w}{k_r}+\frac{w'}{k'_r}-\frac{w w'}{k_rk'_r}\right)^{-1}$$ 
The above equalities follow from ([6] pg. 280 Entry 13-xii) and the definition of $w$. Note that $m$ is the multiplier.\\
Hence for given $0<w<1$ we find $L\in\bf R\rm$ and we get the following parametric evaluation for the Rogers Ramanujan continued fraction
$$R\left(e^{-\pi\sqrt{r(L)}}\right)^{-5}-11-R\left(e^{-\pi\sqrt{r(L)}}\right)^{5}=a_r=$$
\begin{equation}
=\frac{16k_r^6-26k_r^4-wk_r^3+10k_r^2+wk_r}{k_r^4-6k_r^3w-20k_rw^3+15w^2k_r^2-6k_rw+15w^4+w^2}\left(\frac{w}{k_r}+\frac{w'}{k'_r}-\frac{w w'}{k_rk'_r}\right)^{3}
\end{equation}
\[
\]
Thus for a given $w$ we find $L$ and $M$ from (15) and (16). Setting the values of $M$, $L$, $w$ in (14) we get the values of $x$ and $y$ (see Proposition 1). Hence from (24) if we find $k^{(-1)}(x)=r$ we know $R(e^{-\pi\sqrt{r}})$. The clearer result is: 
\[
\]
\textbf{Main Theorem.}\\
When $w$ is a given real number, we can find $x$ from equation (14). Then for the Rogers Ramanujan continued fraction holds
$$R\left(e^{-\pi\sqrt{k^{(-1)}(x)}}\right)^{-5}-11-R\left(e^{-\pi\sqrt{k^{(-1)}(x)}}\right)^{5}=a_r=$$
$$
=\frac{16x^6-26x^4-wx^3+10x^2+wx}{x^4-6x^3w-20xw^3+15w^2x^2-6xw+15w^4+w^2}\times
$$
\begin{equation}
\times\left(\frac{w}{x}+\frac{w'}{\sqrt{1-x^2}}-\frac{w w'}{x\sqrt{1-x^2}}\right)^{3}
\end{equation} 
\[
\]
\textbf{Note.} In the case of $x=k_r$, then $k^{(-1)}(x)=r$ and we have the clasical evaluation with $k_{25r}$ (see [12]). 
\[
\]
\textbf{Theorem 1.} (The first derivative)
$$R'\left(e^{-\pi\sqrt{k^{(-1)}(x)}}\right)=\frac{2^{4/3}x^{1/2}(1-x^2)}{5w^{1/6}w'^{2/3}}
\left(\frac{w}{x}+\frac{w'}{\sqrt{1-x^2}}-\frac{w w'}{x\sqrt{1-x^2}}\right)^{1/2}\times$$ 
\begin{equation}
\times R\left(e^{-\pi\sqrt{k^{(-1)}(x)}}\right)\frac{K^2(x)e^{\pi\sqrt{k^{(-1)}(x)}}}{\pi^2}
\end{equation}
\textbf{Proof.}\\
Combining (8) and (10) and Proposition 2 we get the proof.
\[
\]
We see how the function $k^{(-1)}(x)$ plays the same role in other continued fractions. Here we consider also the Ramanujan's Cubic fraction (see [4]), which is completely solvable using $k_r$.\\
Define the function: 
\begin{equation}
G(x)=\frac{x}{\sqrt{2\sqrt{x}-3x+2x^{3/2}-2\sqrt{x}\sqrt{1-3\sqrt{x}+4x-3x^{3/2}+x^2}}}
\end{equation}
Set for a given $0<w_3<1$ 
\begin{equation}
x=G(w_3)
\end{equation}
Then as in Main Theorem, for the Cubic continued fraction $V(q)$, holds (see [4]):
\begin{equation}
t=V\left(e^{-\pi\sqrt{k^{(-1)}(x)}}\right)=\frac{(1-x^2)^{1/3} w_3^{1/4}}{2^{1/3}x^{1/3}(1-\sqrt{w_3})}
\end{equation}
Observe here that again we only have to know $k^{(-1)}(x)$.\\ 
If happens $x=k_r$, for a certain $r$, then  
\begin{equation}
k_{9r}=\frac{w_3}{k_r}
\end{equation}
and if we set
\begin{equation} 
T=\sqrt{1-8V(q)^3} ,
\end{equation}
then holds
\begin{equation}
(k_r)^2=x^2=\frac{(1-T)(3+T)^3}{(1+T)(3-T)^3}
\end{equation} 
which is solvable always in radicals quartic equation. 
When we know $w_3$ we can find $k_r=x$ from $x=G(w_3)$ and hence $t$.\\
The inverse also holds: If we know $t=V(q)$ we can find $T$ and hence $k_r=x$. The $w_3$ can be find by the degree 3 modular equation which is always solvable in radicals:
$$
\sqrt{k_rk'_r}+\sqrt{k_{9r}k'_{9r}}=1
$$
Let now
\begin{equation}
V(q)=z\Leftrightarrow q=V^{(-1)}(z)
\end{equation}
if
\begin{equation}
V_i(t):=\sqrt{\frac{1-\sqrt{1-8t^3}}{1+\sqrt{1-8t^3}}\left(\frac{3+\sqrt{1-8t^3}}{3-\sqrt{1-8t^3}}\right)^3}
\end{equation}
then 
\begin{equation}
V_i(V(e^{-\pi\sqrt{x}}))=k_x
\end{equation}
or
$$
V(e^{-\pi\sqrt{r}})=V_i^{(-1)}(k_r)
$$
$$
V(e^{-\pi\sqrt{k^{(-1)}(x)}})=V^{(-1)}_i(x)
$$
or
$$e^{-\pi\sqrt{k^{(-1)}(x)}}=V^{(-1)}(V_i^{(-1)}(x))=(V_i\circ V)^{(-1)}(x)$$
and 
\begin{equation}
k^{(-1)}(V_i(V(q)))=\frac{1}{\pi^2}\log(q)^2=r
\end{equation}
Setting now values into (36) we get values for $k^{(-1)}(.)$. The function $V_i(.)$ is an algebraic function.
 
\section{Some Evaluations of the Rogers Ramanujan Continued Fraction}

Note that if $x=k_r$, $r\in\bf Q\rm$ then we have the classical evaluations with $k_r$ and $k_{25r}$.\\
\textbf{Evaluations.}\\
\textbf{1)}
\[
R(e^{-2\pi})=\frac{-1}{2}-\frac{\sqrt{5}}{2}+\sqrt{\frac{5+\sqrt{5}}{2}}
\]
\[
R'(e^{-2\pi})=8\sqrt{\frac{2}{5}\left(9+5\sqrt{5}-2\sqrt{50+22\sqrt{5}}\right)} \frac{e^{2\pi}}{\pi^3}\Gamma\left(\frac{5}{4}\right)^4
\]
\textbf{2)}
Assume that $x=\frac{1}{\sqrt{2}}$, hence   $k^{(-1)}\left(\frac{1}{\sqrt{2}}\right)=1$. From (16) which for this $x$ can be solved in radicals, with respect to $w$, we find  $$w=\frac{\sqrt{2}}{4}\left(\sqrt{5}-1\right)-\frac{1}{2}\sqrt{7\sqrt{5}-15}$$ Hence from $$w'=\sqrt{\sqrt{1-\frac{w^4}{x^2}}\sqrt{1-x^2}}$$ we get $$w'=\left(\frac{1+21\sqrt{-30+14\sqrt{5}}-9\sqrt{-150+70\sqrt{5}}}{\sqrt{2}}\right)^{1/4}$$ Setting these values to (25) we get the value of $a_r$ and then $R(q)$ in radicals. The result is
$$
R(e^{-\pi})^{-5}-11-R(e^{-\pi})^5
=-\frac{1}{8}\left(3+\sqrt{5}-\sqrt{-30+14 \sqrt{5}}\right)[1-\sqrt{5}+
\sqrt{-30+14 \sqrt{5}}+$$
$$+2^{3/8}\left(-3+\sqrt{5}-\sqrt{-30+14 \sqrt{5}}\right) \left(1+21 \sqrt{-30+14 \sqrt{5}}-9 \sqrt{-150+70 \sqrt{5}}\right)^{1/4}]^3\times$$
$$\times[\sqrt{-1574+704\sqrt{5}}-655 \sqrt{-30+14 \sqrt{5}}+293 \sqrt{-150+70 \sqrt{5}}]^{-1}
$$
\textbf{3)}Set $w=1/64$ and $a=1359863889$, $b=36855$, then $$x=$$ 
$$
9 \left(\sqrt{a}+b\right)^{5/6}[49152 6^{1/3} \left(\sqrt{a}+b\right)^{1/6}-960 \left(\sqrt{a}+b\right)^{5/6} +2\cdot6^{2/3}\left(\sqrt{a}+b\right)^{3/2}-
$$ 
$$
-2\cdot6^{5/6}\left(\sqrt{a}+b\right)^{1/6}\surd [-86980957248+36855\cdot 2^{2/3} 3^{1/6}\sqrt{453287963}\cdot(36855+
$$
$$
+\sqrt{a})^{2/3}-2358720 \sqrt{a}+150958080\cdot6^{1/3} \left(\sqrt{a}+b\right)^{1/3}+
$$ 
$$
+4096\cdot 2^{1/3} 3^{5/6} \sqrt{453287963} \left(\sqrt{a}+b\right)^{1/3}+453025819\cdot 6^{2/3} \left(\sqrt{a}+b\right)^{2/3}]+
$$
$$
+384\cdot 2^{2/3} 3^{1/6} \surd[-2358720-64 \sqrt{a}+8192\cdot 6^{1/3} \left(\sqrt{a}+b\right)^{1/3}+
$$
$$
+12285\cdot 6^{2/3} \left(\sqrt{a}+b\right)^{2/3}+2^{2/3} 3^{1/6} \sqrt{453287963}\left(\sqrt{a}+b\right)^{2/3}]^{-1}
$$ 
\textbf{4)} For $$w=\sqrt{\frac{277}{108}+\frac{13\sqrt{385}}{108}}$$ 
we get
$$x=\frac{\sqrt{\frac{277}{12}+\frac{13\sqrt{385}}{12}}}{4+\sqrt{7}}$$
Hence $$R\left(\exp\left[-\pi\cdot k^{(-1)}\left(\frac{\sqrt{\frac{277}{12}+\frac{13\sqrt{385}}{12}}}{4+\sqrt{7}}\right)^{1/2}\right]\right)=$$
$$=\left(-\frac{-8071}{18}+\frac{1075\sqrt{55}}{18}+\frac{1}{18}\sqrt{5(25740148-3470530\sqrt{55})}\right)^{1/5}$$
\textbf{5)} Set $q=e^{-\pi\sqrt{r_0}}$, then from
$$
V(e^{-\pi\sqrt{r_0}})=V^{(-1)}_i(k_{r_0})=V_0
$$ 
and from
$$
V(q^{1/3})=\sqrt[3]{V(q)\frac{1-V(q)+V(q)^2}{1+2V(q)+4V(q)^2}}
$$
We can evaluate all 
$$
V(q_0(n))=b_0(n)=\textrm{Algebraic function of }  r_0
$$ 
where 
$$q_0(n)=e^{-\pi\sqrt{r_0}/3^n}$$
and 
$$V_i(V(q_0(n)))=V_i(b_0(n))=k_{r_0/9^n}$$ 
hence 
$$k^{(-1)}(V_i(b_0(n)))=\frac{r_0}{9^n}$$
An example for $r_0=2$ is $$V(e^{-\pi\sqrt{2}})=-1+\sqrt{\frac{3}{2}}$$  $$V(e^{-\pi\sqrt{2}/3})=\frac{1}{2^{1/3}}\left(-1+\sqrt{\frac{3}{2}}\right)^{1/3}$$
$$V(e^{-\pi\sqrt{2}/9})=\rho_3^{1/3}$$
Where $\rho_3$ can be evaluated in radicals but for simplicity we give the polynomial form.
$$-1-72x-6408x^2+50048x^3+51264x^4-4608x^5+512x^6=0$$
$$\ldots$$
Then respectively we get the values
\begin{equation}
k^{(-1)}\left(-49+35\sqrt{2}+4\sqrt{3(99-70\sqrt{2})}\right)=2/9
\end{equation}
\begin{equation}
k^{(-1)}\left(V_i(\rho_3^{1/3})\right)=2/81
\end{equation}
$$\ldots$$
Hence 
\begin{equation}
k^{(-1)}\left(V_i(b_0(n))\right)=r_0/9^n
\end{equation}   
and  
$$
R\left(e^{-\pi \sqrt{r_0}/3^n}\right)^{-5}-11-R\left(e^{-\pi\sqrt{r_0}/3^n}\right)^{5}=$$
$$
=\frac{16x_n^6-26x_n^4-w_nx_n^3+10x_n^2+w_nx_n}{x_n^4-6x_n^3w_n-20x_nw_n^3+15w_n^2x_n^2-6x_nw_n+15w_n^4+w_n^2}\times
$$
\begin{equation}
\times\left(\frac{w_n}{x_n}+\frac{w'_n}{\sqrt{1-x_n^2}}-\frac{w_n w'_n}{x_n\sqrt{1-x_n^2}}\right)^{3} 
\end{equation}
Where $x_n=V_i(b_0(n))=\textrm{known}$. The $w_n$ are given from (13) (in this case we don't find a way to evaluate $w_n$ in radicals, but as a solution of (13)).\\
6) Set now $$w_0=\frac{-64+a+\sqrt{4096+a (88+a)}}{6 \sqrt{6} \sqrt{a}}$$
then 
$$x_0=\frac{-64+a+\sqrt{4096+a (88+a)}}{\sqrt{6} \left(-4+\sqrt{-2+\frac{16}{a^{1/3}}+a^{1/3}} a^{1/6}+a^{1/3}\right)^2 a^{1/6}}$$
$$
R\left(e^{-\pi\sqrt{k^{(-1)}(x_0)}}\right)^{-5}-11-R\left(e^{-\pi\sqrt{k^{(-1)}(x_0)}}\right)^{5}:=A(a)
$$
where the $A(a)$ is a known algebraic function of $a$ and can  calculated from the Main Theorem. Setting arbitrary real values to $a$ we get algebraic evaluations of the RRCF as in evaluation 4.\\
If we set 
$$
g(x):=\frac{-64+x+\sqrt{4096+x(88+x)}}{\sqrt{6} \left(-4+x^{1/6}\sqrt{-2+\frac{16}{x^{1/3}}+x^{1/3}} +x^{1/3}\right)^2 x^{1/6}}
$$
and if we manage to write $k_r$ in the form $g(a_r)$ for a certain $a_r$ i.e.\\ $V_i(V(e^{-\pi\sqrt{r}}))=k_r=g(a_r)$,
then
$$
R\left(e^{-\pi\sqrt{r}}\right)^{-5}-11-R\left(e^{-\pi\sqrt{r}}\right)^{5}=A(a_r)=A\left(g^{(-1)}(k_r)\right)
$$

\newpage

\centerline{\bf References}\vskip .2in

\noindent

[1]: M.Abramowitz and I.A.Stegun, 'Handbook of Mathematical Functions'. Dover Publications, New York. 1972.

[2]: G.E.Andrews, 'Number Theory'. Dover Publications, New York. 1994.
 
[3]: G.E.Andrews, Amer. Math. Monthly, 86, 89-108(1979).

[4]: Nikos Bagis, 'The complete evaluation of Rogers Ramanujan and other continued fractions with elliptic functions'. arXiv:1008.1304v1.

[5]: Nikos Bagis and M.L. Glasser, 'Integrals related with Rogers Ramanujan continued fraction and q-products'. arXiv:0904.1641. (2009) 

[6]: B.C.Berndt, 'Ramanujan`s Notebooks Part III'. Springer Verlang, New York (1991)

[7]: I.S. Gradshteyn and I.M. Ryzhik, 'Table of Integrals, Series and Products'. Academic Press (1980).

[8]: L. Lorentzen and H. Waadeland, Continued Fractions with Applications. Elsevier Science Publishers B.V., North Holland (1992).  

[9]: S.H.Son, 'Some integrals of theta functions in Ramanujan's lost notebook'. Proc. Canad. No. Thy Assoc. No.5 (R.Gupta and K.S.Williams, eds.), Amer. Math. Soc., Providence.

[10]: E.T.Whittaker and G.N.Watson, 'A course on Modern Analysis'. Cambridge U.P. (1927)

[11]: I.J. Zucker, 'The summation of series of hyperbolic  functions'. SIAM J. Math. Ana.10.192. (1979)

[12]: B.C. Berndt, H.H.Chan, S.S Huang, S.Y.Kang, J.Sohn and S.H.Son, 'The Rogers Ramanujan Continued Fraction'. (page stored in the Web).

\end{document}